\documentclass[a4paper,12pt]{article}

\usepackage{amsmath,amssymb}

\numberwithin{equation}{section}
\newtheorem{thm}{Theorem}[section]

\newtheorem{lem}[thm]{Lemma}
\newtheorem{prp}[thm]{Proposition}

\newtheorem{rmk}[thm]{Remark}
\newtheorem{dfn}[thm]{Definition}
\newtheorem{cnj}[thm]{Conjecture}
\newenvironment{prf}{\noindent {\it Proof} \ }{\hfill $\Box$}
\newenvironment{prfn}[1]{\vskip 1em \noindent {\it Proof of #1} \ }{\hfill $\Box$}

\newcommand{\e}{\epsilon}
\newcommand{\p}{\partial}
\newcommand{\nn}{\nonumber}
\newcommand{\C}{\mathbb{C}}
\newcommand{\F}{\mathcal{F}}
\newcommand{\G}{\mathcal{G}}

\begin{document}

\title{The Inversion Symmetry of the WDVV Equations and Tau Functions}
\author{Si-Qi Liu, Dingdian Xu, Youjin Zhang \\
{\small Department of Mathematical Sciences, Tsinghua University} \\
{\small Beijing 100084, P. R. China}}
\date{\small{\em{Dedicated to Boris Dubrovin on the occasion of his 60th birthday}} }
\maketitle

\begin{abstract}
For two solutions of the WDVV equations that are related by the inversion symmetry,  we show that the associated principal hierarchies of integrable systems
are related by a reciprocal transformation, and the tau functions of the hierarchies are related by a Legendre type transformation. We also consider relationships
between the Virasoro constraints and the topological deformations of the principal hierarchies.
\end{abstract}

\section{Introduction}

The Witten-Dijkgraaf-Verlinde-Verlinde (WDVV) equations, which arise in the study of 2D topological field theory (TFT) in the beginning of 90's of the last century,
are given by the following system of PDEs for an analytic function $F=F(v^1,\dots, v^n)$:
\begin{itemize}
\item[i)] The variable $v^1$ is specified so that
\begin{equation}
\eta_{\alpha\beta}:=\frac{\p^3 F}{\p v^1\p v^\alpha \p v^\beta}=\mbox{constant}, \quad \det(\eta_{\alpha\beta})\ne 0; \label{wdvv-eta}
\end{equation}

\item[ii)] The functions $c^\alpha_{\beta\gamma}:=\eta^{\alpha\nu} c_{\nu\beta\gamma}$ with
\begin{equation}
c_{\alpha\beta\gamma}=\frac{\p^3 F}{\p v^\alpha \p v^\beta\p v^\gamma}, \quad (\eta^{\alpha\beta})=(\eta_{\alpha\beta})^{-1} \label{wdvv-c}
\end{equation}
yield the structure constants of an associative algebra for any fixed $v=(v^1,\dots,v^n)$, i.e, they satisfy
\begin{equation}
c_{\alpha\beta}^\lambda c_{\lambda\gamma}^\nu=c_{\gamma\beta}^\lambda c_{\lambda\alpha}^\nu, \quad \mbox{for any } 1\le \alpha, \beta, \gamma, \nu \le n. \label{wdvv}
\end{equation}
Here and in what follows summation with respect to repeated upper and lower indices is assumed.
\end{itemize}
In \cite{D3, D1} Dubrovin formulated the WDVV equations with an additional quasi-homogeneity condition on $F$, we will recall this condition in Sec. \ref{sec-4}
and call a solution of \eqref{wdvv-eta}--\eqref{wdvv} satisfying the quasi-homogeneous condition a conformal solution of the WDVV equations. These equations
together with the quasi-homogeneity condition are satisfied by the primary free energy $F$ of the matter sector of a 2D TFT with $n$
primary fields as a function of the coupling constants \cite{Dij2, Dij3, Witten1}. In \cite{D3, D1} Dubrovin reformulated these equations in a coordinate
free form by introducing the notion of Frobenius manifold structure on the space of the parameters $v^1,\dots, v^n$ and revealed significantly rich geometric structures
of the WDVV equations, which have become important in the study of several different areas of mathematical research, including the theory of Gromov - Witten
invariants, singularity theory and nonlinear integrable systems, see \cite{D1, weyl, DZ, DZ2} and references therein. In particular, such geometrical structures enable one to associate
a solution of the WDVV equations with a hierarchy of bihamiltonian integrable systems of hydrodynamic type which is called the principal hierarchy in
\cite{DZ2}. This hierarchy of integrable systems plays important role in the procedure of reconstructing a 2D TFT from its primary
free energy as a solution of the WDVV equations. In this construction, the tau function that corresponds to a particular solution of the principal hierarchy
serves as the genus zero partition function, and the full genera partition function of the 2D TFT is a particular tau function of an integrable hierarchy of
evolutionary PDEs of KdV type which is certain deformation of the principal hierarchy, such a deformation of the principal hierarchy is call the
{\em{topological deformation}} \cite{DZ2}.

In this paper we are to interpret certain discrete symmetry of the WDVV equations in term of the associated principal hierarchy and its tau function. The
discrete symmetry we consider here was given by Dubrovin in  Appendix B of \cite{D1} and is called the inversion symmetry. This symmetry is obtained
from a special Schlesinger transformation of the system of linear ODEs with rational coefficients associated to the WDVV equations (see Remark 4.2 of \cite{D2} for details). It turns out that in terms of the principal hierarchy and its tau function the inversion symmetry has a
simple interpretation. On the principal hierarchy it acts as certain reciprocal transformation, and on the associated tau function it acts as a Legendre type
transformation. In Appendix B of \cite{D1} there is given another class of discrete symmetries of the WDVV equations which are called the Legendre-type
transformations or the type I symmetries, the relation of such symmetries with the principal hierarchy and its tau function is given in \cite{DZ2}. Besides these discrete symmetries,
the WDVV equations also possess continuous symmetries whose Lie algebra of infinitesimal generators (without the quasi-homogeneity condition) was studied in
\cite{CKS, FLS, Givental, Lee, Leur}.

Recall \cite{D1} that  a symmetry of the WDVV equations is given by a transformation
\begin{equation}
v^\alpha\mapsto \hat{v}^\alpha,\quad
\eta_{\alpha\beta}\mapsto\hat{\eta}_{\alpha\beta},\quad F\mapsto \hat{F}
\end{equation}
that preserves the WDVV equations. The inversion symmetry
given in \cite{D1} has the following form:
\begin{align}
\begin{split}
&\hat{v}^1=\frac12\,\frac{\eta_{\alpha\beta}  v^\alpha v^\beta}{v^n},\quad \hat{v}^i=\frac{v^i}{v^n}\ (i=2,\dots, n-1),\quad \hat{v}^n=-\frac1{v^n},\\
&\hat{\eta}_{\alpha\beta}=\eta_{\alpha\beta},\quad \hat{F}(\hat{v})=(v^n)^{-2}\left(F(v)-\frac12 \eta_{\alpha\beta}v^1v^\alpha v^\beta\right).
\end{split}\label{IS}
\end{align}
Here we assume that $\eta_{11}=0$ and the coordinates $v^1,\dots, v^n$ are normalized such that the constants $\eta_{\alpha\beta}$ take the values
\begin{equation}
\eta_{\alpha\beta}=\delta_{\alpha+\beta,n+1}.
\end{equation}
This can always be achieved by performing an invertible linear transformation
\[v^1\mapsto \tilde{v}^1=v^1+\sum_{i=2}^n b_i v^i,\quad v^j \mapsto \tilde{v}^j=\sum_{i=2}^n a^j_i v^i\ (j=2, \dots, n),\]
where $b_i, a^j_i \in \C$.
We call the solution $\hat{F}(\hat{v})$ of the WDVV equations \eqref{wdvv-eta}--\eqref{wdvv} the inversion of the solution $F(v)$.

The content of the paper is arranged as follows. We first recall in Sec.\,\ref{sec-2} the definition of the principal hierarchy and its tau functions associated to
a calibrated solution of the WDVV equations. We then show in Sec.\,\ref{sec-3} that the action of the inversion symmetry of the WDVV equations on principal
hierarchies is given by certain reciprocal transformation, and we give the transformation rule of the associated tau functions, see Propositions \ref{prp-a}
and \ref{prp-b}. In Sec.\,\ref{sec-4} we impose the conformal condition on solutions of the WDVV equations and consider the transformation rule of the
inversion symmetry on principal hierarchies and their bihamiltonian structures. These results will be used when we study the topological deformations of
principal hierarchies. In Sec.\,\ref{sec-5} we consider transformation rule of the Virasoro constraints for tau functions of the principal hierarchies. In Sec.\,\ref{sec-6} we consider the action of the inversion symmetry on topological deformations of the principal hierarchies and their tau functions.

The present paper is a rewritten of an early preprint \cite{XZ}. We omit the content associated to the type I symmetries of the WDVV equations, and refine
presentations of the results given there. The main new content is the proof of Conjecture 6.1 of \cite{XZ}.

\section{Calibrations, Principal Hierarchies, and \\ Tau Functions}\label{sec-2}

The notion of calibrations of a solution of the WDVV equations (or a Frobenius manifold) corresponds to the choice of a system of deformed flat coordinates on a Frobenius manifold \cite{D1}, it was first introduced in \cite{Givental-a} and then modified in \cite{CKS}. In what follows we use the modified one.

It is well-known that the system of WDVV equations is equivalent to the flatness of the deformed flat connection (see \cite{D1} for details):
\[\tilde{\nabla}_X(z)Y=\nabla_XY+z\,X\cdot Y,\]
where $X, Y$ are vector fields on the Frobenius manifold $M$, $\nabla$ is the Levi-Civita connection of the metric $(\eta_{ij})$,
and $z$ is an arbitrary nonzero complex number.

The flatness of the deformed connection implies the existence of deformed flat coordinates, so the following equation
\begin{equation}
\tilde{\nabla}d\tilde{v}(v,z)=0 \label{flat-coord}
\end{equation}
has $n$ linearly independent solutions $\tilde{v}^1(v,z), \dots, \tilde{v}^n(v,z)$ 
which are analytic at $z=0$.  We denote them by
\[\tilde{v}^{\alpha}(v,z)=\eta^{\alpha\nu}\theta_{\nu}(z)=\eta^{\alpha\nu}\sum_{p\ge0}\theta_{\nu,p}(v)z^p,\]
then the equation \eqref{flat-coord} implies that
\begin{equation}
\p_{\alpha}\p_{\beta}\theta_{\nu}(z)=z\,c^{\gamma}_{\alpha\beta}\p_{\gamma}\theta_{\nu}(z),\  \ \p_{\alpha}=\frac{\p}{\p v^{\alpha}},\  \alpha, \beta, \nu=1, \dots, n. \label{eq-theta}
\end{equation}
Solutions to the above equations are not unique,  most of the results given below in this section hold true only for those solutions that are normalized by certain conditions coming from topological field theories. These carefully chosen solutions are called calibrations of the Frobenius manifolds $M$.

\begin{dfn}
Let $F(v)$ be a solution of the WDVV equations, a family of functions
\[\{\theta_{\alpha, p}(v)\mid \alpha=1, \dots, n;\ p=0, 1, 2, \dots\}\]
is called a calibration of $F(v)$ if their generating functions
\[\theta_{\alpha}(z)=\sum_{p\ge0} \theta_{\alpha, p}(v)z^p\]
satisfy the equations \eqref{eq-theta} and the normalization conditions
\begin{align}
&\theta_{\alpha}(0)=v_\alpha:=\eta_{\alpha\beta}v^{\beta},\label{zh-10-5} \\
&\p_{\mu}\theta_{\alpha}(z)\,\eta^{\mu\nu}\,\p_{\nu}\theta_{\beta}(-z)=\eta_{\alpha\beta},\\
&\theta_{\alpha,1}(v)=\frac{\p F}{\p v^\alpha},\\
&\p_1 \theta_{\alpha}(z)=z\,\theta_{\alpha}(z)+\eta_{1\alpha}.
\end{align}
The solution $F(v)$ together with a calibration $\{\theta_{\alpha,p}(v)\}$ is called a calibrated solution of the WDVV equations.
\end{dfn}

Let $(F(v), \{\theta_{\alpha,p}(v)\})$ be a calibrated solution of the WDVV equations, we introduce a hierarchy of evolutionary PDEs of  hydrodynamic type:
\begin{equation}
\frac{\p v^\gamma}{\p t^{\alpha, p}}=\eta^{\gamma\beta}\frac{\p}{\p x}\left(\frac{\p \theta_{\alpha, p+1}}{\p v^{\beta}}\right),
\quad \alpha, \gamma=1, \dots, n,\ p\ge0. \label{PH}
\end{equation}
It is easy to see that
\[\frac{\p v^\gamma}{\p t^{1,0}}=\frac{\p v^\gamma}{\p x},\]
so in what follows we identify $t^{1,0}$ with $x$. By using the WDVV equations one can prove the following results:
\begin{itemize}
\item[i)] Each flow  $\p_{\alpha,p}$ possesses a Hamiltonian formalism with the Hamiltonian operator $P_1=(P_1^{\alpha\beta})=\eta^{\alpha\beta}\p_x$ and Hamiltonian
$H_{\alpha, p}[v]=\int \theta_{\alpha, p+1}(v)\,dx$;
\item[ii)] $\{H_{\alpha,p}, H_{\beta,q}\}_1=0$, where $\{\ ,\ \}_1$ is the Poisson bracket defined by 
\begin{equation}
\{H_1, H_2\}_1=\int \frac{\delta H_1}{\delta v^\alpha}\,P^{\alpha\beta}_1\left(\frac{\delta H_2}{\delta v^\beta}\right)dx; \label{poi-bra}
\end{equation}
\item[iii)] Denote  $\p_{\alpha,p}:=\frac{\p}{\p t^{\alpha,p}}$, then $\p_{{\beta,q}}\theta_{\alpha, p}(v)=\p_{\alpha, p}\theta_{\beta,q}(v).$
\end{itemize}
The second property ii) also implies
\[[\p_{\alpha,p}, \p_{\beta,q}]=0,\quad \p_{\alpha,p}H_{\beta, q}=0,\]
which are easy corollaries of the properties of Poisson brackets.

\begin{dfn}
The hierarchy \eqref{PH} of integrable evolutionary PDEs is called the principal hierarchy associated to the calibration $\{\theta_{\alpha,p}(v)\}$.
\end{dfn}
\begin{rmk} The notion of principal hierarchy was introduced in \cite{DZ2} for Frobenius manifolds associated to conformal solutions of the WDVV
equations, in this case calibrations are given by the deformed flat coordinates of the
Frobenius manifolds \cite{D1}, see Sec. \ref{sec-4} below.
\end{rmk}

Now we are to define tau functions of the principal hierarchy as it is done in \cite{D1, DZ2}. First we recall the definition \cite{D1} of the family of functions
\[\{\Omega_{\alpha, p;\beta,q}(v)\mid \alpha,\beta=1, \dots, n;\ p,q=0, 1, 2, \dots\}\]
by the following generating functions
\begin{equation}\label{omega}
\sum_{p,q\ge0}\Omega_{\alpha, p;\beta,q}(v)\,z^p\,w^q=\frac{\p_{\mu}\theta_{\alpha}(z)\,\eta^{\mu\nu}\,\p_{\nu}\theta_{\beta}(w)-\eta_{\alpha\beta}}{z+w}.
\end{equation}
Then one can prove that
\begin{equation}
\Omega_{\alpha, p;\beta,q}=\Omega_{\beta,q;\alpha, p}, \quad \p_{\gamma,s}\Omega_{\alpha, p;\beta,q}=\p_{\alpha, p}\Omega_{\gamma,s;\beta,q}, \label{omega-cond}
\end{equation}
which imply that if $v^\alpha(t)$ is a solution of the principal hierarchy associated to certain calibration
(here we use $t$ to denote all the time variables $t^{\alpha,p}$ of the principal hierarchy),
then there exists a function $f(t)$ such that
\[\Omega_{\alpha, p;\beta,q}(v(t))=\p_{\alpha, p}\p_{\beta,q}f(t).\]
In particular, we have $\theta_{\alpha,p}(v(t))=\Omega_{1,0;\alpha, p}(v(t))=\p_{1,0}\p_{\alpha, p}f(t)$.

\begin{dfn}\label{dfn-tau}
Let $(F(v), \{\theta_{\alpha,p}(v)\})$ be a calibrated solution of the WDVV equations, $\{\p_{\alpha,p}\}$ be the associated principal hierarchy.
A function $\tau(t)$ is called a tau function of the principal hierarchy if
\begin{equation}\label{tau-def}
\Omega_{\alpha, p;\beta,q}(v(t))=\p_{\alpha, p}\p_{\beta,q}\log \tau(t),
\end{equation}
where $v^\alpha(t)=\eta^{\alpha\beta}\p_{1,0}\p_{\beta,0}\log \tau(t)$.
\end{dfn}

Note that if $\tau(t)$ is a tau function of the principal hierarchy, then $v^\alpha(t)=\eta^{\alpha\beta}\p_{1,0}\p_{\beta,0}\log \tau(t)$ is a solution of the principal
hierarchy. Indeed, by using the property of the functions $\Omega_{\alpha,p;\beta,q}$ \cite{D1}
\begin{equation}\label{zh-22}
\p_\xi\Omega_{\alpha,p;\beta,q}=\frac{\p\theta_{\alpha,p}}{\p v^\sigma}\frac{\p\theta_{\beta,q}}{\p v^\lambda} c^{\sigma\lambda}_\xi,\quad
c^{\sigma\lambda}_\xi=\eta^{\sigma\gamma} c^\lambda_{\gamma\xi}
\end{equation}
we have
\begin{eqnarray}
\frac{\p v^\alpha(t)}{\p t^{\beta,q}}&=&\eta^{\alpha\gamma}
\p_{1,0}\p_{\gamma,0}\p_{\beta,q}\log \tau(t)
=\eta^{\alpha\gamma}
\p_{1,0}\Omega_{\gamma,0;\beta,q}(v(t))\nn\\
& =&\eta^{\alpha\gamma} \frac{\p\Omega_{\gamma,0;\beta,q}}{\p v^\xi}
v^\xi_x(t)=\eta^{\alpha\gamma} \frac{\p\theta_{\beta,q}}{\p v^\lambda} c^\lambda_{\gamma\xi} v^\xi_x(t)\nn\\
&=&\eta^{\alpha\gamma} \p_x  \frac{\p\theta_{\beta,q+1}}{\p v^\gamma}.\nn
\end{eqnarray}
On the other hand, the above argument shows that a solution of the principal hierarchy
also defines a tau function.

\section{The Actions of the Inversion Symmetry}\label{sec-3}

In this section, we study the actions of the inversion symmetry on calibrations, principal hierarchies, and tau functions of a solution of the WDVV equations.
We fix a pair of solutions $F(v), \hat{F}(\hat{v})$ of the WDVV equations that are related by the transformation \eqref{IS}.

\begin{prp}
Let $\{\theta_{\alpha,p}(v)\}$ be a calibration of $F(v)$, then the following functions
\begin{align}\label{ntheta}
&\hat{\theta}_{1,0}(\hat{v})=-\frac{1}{v^n},\quad \hat{\theta}_{1,p}(\hat{v})=-\frac{\theta_{n,p-1}(v)}{v^n},\quad  p\ge1,\nn \\
&\hat{\theta}_{i,p}(\hat{v})=\frac{\theta_{i,p}(v)}{v^n},\quad
 2\le i\le n-1,\ p\ge 0, \\
&\hat{\theta}_{n,p}(\hat{v})=\frac{\theta_{1,p+1}(v)}{v^n},\quad p\ge 0,\nn
\end{align}
give a calibration $\{\hat{\theta}_{\alpha,p}(\hat{v})\}$ of $\hat{F}(\hat{v})$.
\end{prp}
\begin{prf}
According to the definition of calibration, we need to show that the generating functions $\hat{\theta}_\alpha(z)$ satisfy the
equation \eqref{eq-theta}, i.e. 
\[\frac{\p^2 \hat{\theta}_\nu(z)}{\p\,\hat{v}^\alpha\p\,\hat{v}^\beta}=z\,\hat{c}_{\alpha\beta}^\gamma
\frac{\p \hat{\theta}_\nu(z)}{\p\,\hat{v}^\gamma},\]
and the normalization conditions. Note that both $\hat{v}$ and $\hat{\theta}_{\alpha,p}$ are defined in three cases,
so to prove the above identity one needs to verify $3\times2\times3=18$ cases. The computation is not
hard (in fact, there is nothing more than elementary calculus), but very lengthy, so we omit the details.
 
The following identity is frequently used in these verifications:
\[v^n\delta_{\alpha}^n\frac{\p v^{\mu}}{\p \hat{v}^\beta}+v^n\delta^n_\beta\frac{\p v^\mu}{\p \hat{v}^\alpha}
=\frac{\p^2 v^\mu}{\p \hat{v}^\alpha\p \hat{v}^\beta}+\eta_{\alpha\beta}v^n\delta_1^\mu,\]
it can be proved by definition and case by case verifications (18 case again). 
By using the above identity and the following ones
\[\hat{c}_{\alpha\beta\gamma}(\hat{v})=\frac{\p^3 \hat{F}(\hat{v})}{\p\hat{v}^\alpha\p\hat{v}^\beta\hat{v}^\gamma}
=(v^n)^{-2}\frac{\p v^\lambda}{\p \hat{v}^\alpha}\frac{\p v^\mu}{\p \hat{v}^\beta}\frac{\p v^\nu}{\p \hat{v}^\gamma}\,c_{\lambda\mu\nu}(v)\]
which was given in \cite{D1}, one can prove the proposition straightforwardly.
\end{prf}

\begin{prp}\label{prp-a}
Let $\{\p_{\alpha,p}\}$ be the principal hierarchy associated to a calibration $\{\theta_{\alpha,p}(v)\}$ of $F(v)$.
Introduce the following reciprocal transformation
\begin{align}
&d\hat{x}=\sum_{\alpha=1}^n\sum_{p\ge0}\theta_{\alpha,p}(v)dt^{\alpha,p},\label{zh-24-1}\\
\begin{split}\label{zh-24-2}
&\hat t^{1,0}=\hat{x},\quad \hat{t}^{1, p}=-t^{n,p-1},\quad p \ge 1,\\
&\hat{t}^{n,p}= t^{1,p+1}\ (p \ge 0),\quad \hat{t}^{i,p}=t^{i,p},\quad 2\le i\le n,\ p \ge 0,
\end{split}
\end{align}
and denote $\hat{\p}_{\alpha,p}:=\frac{\p}{\p \hat{t}^{\alpha,p}}$, then we have
\begin{equation}\label{NPH}
\frac{\p\hat{v}^\beta}{\p \hat{t}^{\alpha,p}}=\hat{\eta}^{\beta\gamma}\frac{\p}{\p\hat{x}}
\left(\frac{\p\hat{\theta}_{\alpha,p+1}(\hat{v})}{\p\hat{v}^\gamma}\right),\ \alpha,\beta=1,\dots,n,\ p\ge 0,
\end{equation}
i.e. $\{\hat{\p}_{\alpha,p}\}$ is the principal hierarchy of $\hat{F}(\hat v)$ with calibration $\{\hat{\theta}_{\alpha,p}(\hat v)\}$.
\end{prp}

\begin{prf}
From the definition of the reciprocal transformation we have
\begin{align}
\begin{split}\label{zh-19}
\frac{\p}{\p\hat x}&=\frac{1}{v^n}\frac{\p}{\p x},\\
\frac{\p}{\p\hat{t}^{1,p}}&=-\frac{\p}{\p t^{n,p-1}}+\frac{\theta_{n,p-1}(v)}{v^n}\frac{\p}{\p x},\quad  p\ge 1,\\
\frac{\p}{\p\hat{t}^{i,p}}&=\frac{\p}{\p t^{i,p}}-\frac{\theta_{i,p}(v)}{v^n}\frac{\p}{\p x},\quad  2\le i\le n,\ p\ge 0,\\
\frac{\p}{\p\hat{t}^{n,p}}&=\frac{\p}{\p t^{1,p+1}}-\frac{\theta_{1,p+1}(v)}{v^n}\frac{\p}{\p x},\quad  p\ge 0.
\end{split}
\end{align}
The proposition can be proved by direct calculation.
\end{prf}

\begin{prp}\label{prp-b}
Let $\tau(t)$ be a tau function of the principal hierarchy associated to a calibration $\{\theta_{\alpha,p}(v)\}$ of $F(v)$. Define
\begin{equation}\label{legendre}
\log \hat{\tau}(\hat{t})=\log \tau(t)-x\frac{\p \log \tau(t)}{\p x},
\end{equation}
then $\hat{\tau}(\hat{t})$ is a tau function of the principal hierarchy associated to $\{\hat{\theta}_{\alpha,p}(\hat{v})\}$.
\end{prp}
\begin{prf}
By definition the functions $\hat{\Omega}_{\alpha,p;\beta,q}(\hat{v})$ are given by \eqref{omega} in terms of the functions $\{\hat{\theta}_{\alpha,p}(\hat{v})\}$.
From  the relation \eqref{IS}, \eqref{ntheta} it follows that
\begin{align}
\hat{\Omega}_{\alpha,p;\beta,q}(\hat{v})=&(-1)^{\delta^1_\alpha+\delta^1_\beta}
\left(\Omega_{\alpha+(n-1)\delta(\alpha),p-\delta(\alpha);\beta+(n-1)\delta(\beta),q-\delta(\beta)}(v)\right.\nn\\
&\quad -\left.\frac{1}{v^n}\, \theta_{\alpha+(n-1)\delta(\alpha),p-\delta(\alpha)}(v)\,\theta_{\beta+(n-1)\delta(\beta),q-\delta(\beta)}(v)\right),
\label{relation}
\end{align}
where $\delta(\alpha)=\delta^1_\alpha-\delta^n_\alpha$, and we assume that $\theta_{n,-1}=1$, $\Omega_{\alpha,p;n,-1}=\Omega_{n,-1;\beta,q}=0$ when
$(\alpha,p), (\beta,q)\ne (1,0)$. Then one can verify, by using \eqref{IS},
\eqref{ntheta} and \eqref{tau-def}, that
\[\hat{\Omega}_{\alpha,p;\beta,q}(\hat{v}(\hat{t}))=\hat{\p}_{\alpha,p}\hat{\p}_{\beta,q}\log \hat{\tau}(\hat{t}),\]
where $\hat{v}^\alpha(\hat{t})=\eta^{\alpha\beta}\hat{\p}_{1,0}\hat{\p}_{\beta,0}\log \hat{\tau}(\hat{t})$. The proposition is proved.
\end{prf}

Let $\tau(t)$ be a tau function of a principal hierarchy, then the reciprocal transformation \eqref{zh-24-1} can be written as
\[d\hat{x}=d\left(\frac{\partial}{\partial x}\log\tau(t)\right).\]
It follows that up to the addition of a constant we have
\begin{equation}
\hat{x}=\frac{\partial}{\partial x}\log\tau(t). \label{zh-20b}
\end{equation}
The constant can be absorbed by a translation of $\hat{x}$ in the definition of the reciprocal transformation, so we will assume from now on the validity of
\eqref{zh-20b}. Thus in terms of a given tau function, the reciprocal transformation \eqref{zh-24-1}, \eqref{zh-24-2} can be represented by \eqref{zh-24-2},
\eqref{legendre} and \eqref{zh-20b}.

We note that the inverse of the transformation \eqref{zh-24-2}, \eqref{zh-20b}, \eqref{legendre} is given by \eqref{zh-24-2} and
\begin{equation}\label{zh-10-10}
x=-\frac{\p}{\p\hat{x}}\log\hat{\tau}(\hat{t}),\quad \log\tau(t)=\log\hat{\tau}(\hat{t})-\hat{x} \frac{\p\log\hat{\tau}(\hat{t})}{\p\hat{x}}.
\end{equation}
They are transformations of Legendre type.

\section{Conformal Case}\label{sec-4}
In this section we are to include the quasi-homogeneity condition into the WDVV equations as it is formulated in \cite{D1}.

\begin{dfn}
A solution $F(v)$ of the WDVV equations is called conformal if there exists a vector field $E$, called the Euler vector field, of the form
\[E=\sum_{\alpha=1}^n\left(q^\alpha_\beta v^\beta+r^\alpha\right)\p_\alpha,\quad q^\alpha_\beta, r^\alpha\in\C,\]
and some constants $d, A_{\alpha\beta}, B_{\alpha}, C\in\C$ such that
\[E(F)=(3-d)F+\frac12 A_{\alpha\beta}v^{\alpha}v^{\beta}+B_{\alpha} v^{\alpha}+C.\]
\end{dfn}

It is often assumed that the matrix $Q=(q^\alpha_\beta)$ is diagnolizable and $q^1_1=1$. The coordinates $v^1,\dots, v^n$ are normalized so that
\[E=\sum_{\alpha=1}^n\left(d_\alpha v^\alpha+r^\alpha\right)\p_\alpha,\quad d_1=1,\]
and $r^\alpha=0$ if $d_\alpha \ne 0$. In this paper, we assume that
\[r^1=\dots=r^n=0.\]
This assumption ensures that the solution $\hat F(\hat v)$ of the WDVV equations
obtained from $F(v)$ by the action of the inversion symmetry also has a diagnolizable
Euler vector field, while this is not always true without the above assumption, see \cite{D1} and Lemma \ref{symm-lm}.
Then the Euler vector field can be written in the following form:
\[E=\sum_{\alpha=1}^n\left(1-\frac{d}2-\mu_{\alpha}\right)v^\alpha\,\p_\alpha,\quad \mu_1=-\frac{d}2,\]
where the constants $d$ and $\{\mu_{\alpha}\}$ are called the charge and the spectum of $F(v)$ respectively \cite{D1}.

Note that the WDVV equations only involve the third order derivatives of $F(v)$, so we can add certain quadratic functions of $v^1,\dots, v^n$ to $F(v)$ such that
the constants $A_{\alpha\beta}, B_{\alpha}, C$ satisfy the following normalizing conditions
\begin{align*}
&A_{\alpha\beta}\ne 0 \mbox{ only if } \mu_\alpha+\mu_\beta=-1,\\
&B_\alpha\ne 0 \mbox{ only if } \mu_\alpha=\frac{d}2-2,\\
&C\ne 0 \mbox{ only if } d=3.
\end{align*}
Further more, our assumption on $\eta_{11}$ and $r^{\alpha}$ implies that
\[A_{1\alpha}=0,\quad B_1=0.\]

The following lemma is proved in \cite{D1}.
\begin{lem}[\cite{D1}]\label{symm-lm}
Let $F(v)$ be a conformal solution of the WDVV equations with charge $d$ and spectrum $\{\mu_\alpha\}$, and $\hat{F}(\hat{v})$ be its inversion,
then $\hat{F}(\hat{v})$ is also conformal, whose charge $\hat{d}$ and spectrum $\{\hat{\mu}_\alpha\}$ read
\begin{equation}
\hat{d}=2-d,\ \hat\mu_1=\mu_n-1,\ \hat{\mu}_n=\mu_1+1,\ \hat{\mu}_i=\mu_i\ (2 \le i \le n-1).
\end{equation}
\end{lem}

\begin{dfn}
Let $F(v)$ be a conformal solution of the WDVV equations with spectrum $\{\mu_\alpha\}$, a calibration $\{\theta_{\alpha,p}(v)\}$ is called conformal
if there exist constant matrices $R_1, R_2, \dots$ such that
\begin{equation}
E\left(\p_\beta\theta_{\alpha,p}(v)\right)=\left(p +\mu_\alpha+\mu_\beta\right)\p_\beta\theta_{\alpha,p}(v)+
\sum_{k=1}^p\p_\beta\theta_{\gamma,p-k}(v)\,\left(R_k\right)^\gamma_\alpha,
\end{equation}
and
\begin{align}
&(R_k)^\alpha_\beta\ne 0 \mbox{ only if } \mu_\alpha-\mu_\beta=k,\label{p-1}\\
&\eta_{\alpha\gamma} (R_k)^\gamma_\beta+(-1)^k \eta_{\beta\gamma}(R_k)^\gamma_\alpha=0.\label{p-2}
\end{align}
\end{dfn}

The property \eqref{p-1} implies that there is only finitely many nonzero matrices $R_k$. These matrices and the metric $(\eta_{\alpha\beta})$,
the spectrum $\{\mu_{\alpha}\}$ form a representative of the monodromy data of $F(v)$ at $z=0$, and
\[v_\alpha(z)=\sum_{p\ge 0} \theta_{\gamma,p}(v) \left(z^\mu z^R\right)^{\gamma}_\alpha,
\quad \alpha=1,\dots, n\]
form a system of flat coordinates for the deformed flat connection of the
Frobenius manifold associated to the conformal solution $F$ of the WDVV equations.
Here $R=R_1+R_2+\dots$ and $\mu=\textrm{diag}(\mu_1,\dots,\mu_n)$. See \cite{D1, D2} for details.

\begin{prp}
Let $(F(v), \{\theta_{\alpha,p}(v)\})$ be a calibrated conformal solution of the WDVV equations, then the calibration $\{\hat{\theta}_{\alpha,p}(\hat{v})\}$
of $\hat{F}(\hat{v})$ is also conformal.
\end{prp}
\begin{prf}
We only need to compute $\hat{E}\left(\hat{\p}_\beta\hat{\theta}_{\alpha,p}(\hat{v})\right)$, then the matrices $\hat{R}_1, \hat{R}_2, \dots$ for $\hat{F}(\hat{v})$
can be obtained
\[(\hat{R}_{k})^{\alpha}_{\beta}=(-1)^{\delta^\alpha_1+\delta^1_\beta}
\left(R_{k+\delta(\alpha)-\delta(\beta))}\right)^{\alpha+(n-1)\delta(\alpha)}_{\beta+(n-1)\delta(\beta)},\quad k=1,2, \dots.\]
The proposition is proved.
\end{prf}

The principal hierarchy associated to a conformal calibtation has a very important additional structure -- the bihamiltonian structure.
We already know from Sec.\,\ref{sec-2} that the principal hierarchy has one Hamiltonian structure $P_1$. When the calibrated solution
$(F(v), \{\theta_{\alpha,p}(v)\})$ is conformal, we have the following results.
\begin{lem}[\cite{D1}]
Define a matrix differential operator $P_2=(P_2^{\alpha\beta})$, where
\begin{align}
&P_2^{\alpha\beta}=g^{\alpha\beta}(v)\p_x+\Gamma^{\alpha\beta}_{\gamma}(v)\,v^{\gamma}_x, \label{p2-1}\\
&g^{\alpha\beta}(v)=\left(1-\frac{d}2-\mu_{\gamma}\right)v^{\gamma}c^{\alpha\beta}_{\gamma}(v),
\ \Gamma^{\alpha\beta}_{\gamma}(v)=\left(\frac12-\mu_{\beta}\right)c^{\alpha\beta}_{\gamma}(v), \label{p2-2}
\end{align}
then $P_2$ is a Hamiltonian operator which is compatible with $P_1$. Furthermore, for any conformal calibration $\{\theta_{\alpha,p}(v)\}$ of $F(v)$, we have
\begin{equation}
\{\cdot\,, H_{\alpha,p-1}\}_2=\left(p+\mu_{\alpha}+\frac12\right)\{\cdot\,, H_{\alpha,p}\}_1+\sum_{k=0}^p \left(R_k\right)^\beta_\alpha\{\cdot\,, H_{\beta,p-k}\}_1,
\end{equation}
where $\{\ ,\ \}_2$ is the Poisson bracket defined by $P_2$, see \eqref{poi-bra}.
\end{lem}

It has been shown in Proposition \ref{prp-a} that the inversion symmetry preserves the first Hamiltonian structure $P_1$, then it is natural to ask: does
it also preserve the second Hamiltonian structure $P_2$?

\begin{prp} \label{prp-reci}
Let $F(v),\ \hat{F}(\hat{v})$ be a pair of solutions of the WDVV equations that are related by the inversion symmetry, and $P_i,\ \hat{P}_i\ (i=1,2)$ be the
corresponding Hamiltonian structures. Denote by $\Phi$ the reciprocal transformation \eqref{zh-24-1} and \eqref{zh-24-2}, then we have
\[\Phi(P_1)=\hat{P}_1,\quad \Phi(P_2)=\hat{P}_2.\]
\end{prp}

We note that
the action of reciprocal transformations of the form \eqref{zh-24-1} and \eqref{zh-24-2} on evolutionary PDEs of hydrodynamic type and their Hamiltonian structures of the form
\eqref{p2-1}
was first investigated by Ferapontov and Pavlov in \cite{Fer}.  After the action of a reciprocal transformation a Hamiltonian operator of the form \eqref{p2-1}
becomes nonlocal in general, the nonlocal Hamiltonian operator is given by a differential operator of the form \eqref{p2-1} plus an integral operator, in this case the metric $(g^{\alpha\beta})$ is no longer flat.
In \cite{abenda} Abenda considered the conditions under which such a reciprocal transformation preserves the locality of a Hamiltonian structures of hydrodynamic type.
In \cite{SZ} we studied a general class of nonlocal Hamiltonian structures in terms of infinite dimensional Jacobi structures and gave the transformation rule of such Hamiltonian structures under certain reciprocal transformations, a criterion on whether
a reciprocal transformation preserves the locality of a Hamiltonian structure was also given in \cite{SZ}.

\begin{thm}[\cite{SZ}]\label{locality}
Let $P$ be a quasi-local bivector, $\rho$ be an invertible differential polynomial, the reciprocal transformation defined by $\rho$
is denoted by $\Phi$. Then $\Phi(P)$ is local if and only if
\[[P,\ \Lambda]=0, \quad z(P, \rho)=0,\]
where $\Lambda=\int \rho\,dx$, and $z(P, \rho)$ is the nonlocal charge of the pair $(P, \rho)$. 
\end{thm}

In the above theorem the bracket $[\,,\,]$ is the Schouten-Nijenhuis bracket defined on the space of quasi-local multi-vectors \cite{SZ}, and the function $\rho$ in our present case is given by $\rho=v^n$. The definition of the nonlocal charge $z(P,\rho)$ will be given below in the proof of Proposition \ref{prp-reci}.

\begin{prfn}{Proposition \ref{prp-reci}}
According to the general results of \cite{SZ} (c.f. \cite{Fer}), $\Phi(P_i)\ (i=1, 2)$ are Jacobi structures of hydrodynamic type. To prove the proposition,
one need to show that they are both local, and their associated metrics coincide with the ones of $\hat{P}_i\ (i=1,2)$.

We first give the proof of the locality of $\Phi(P_2)$ by using Theorem \ref{locality}.
The proof of locality for $\Phi(P_1)$ is easier and is omited here.

In our reciprocal transformation \eqref{zh-24-2}, $\rho=v^n$, so $\Lambda=\int v^n\,dx$,
then we need to show that
\begin{equation}\label{zh-10-1}
[P_2, \Lambda]=0,
\end{equation}
which is equivalent to say that there exists a constant $c$ such that
\[\nabla^i\nabla_k\left(v^n\right)=c\,\delta^i_k,\]
where $\nabla$ is the Levi-Civita connection of the metric $(g_{\alpha\beta})=(g^{\alpha\beta})^{-1}$ (see \eqref{p2-2}).
By a straightforward calculation one can obtain that $c=\frac{1-d}2$, thus \eqref{zh-10-1}
holds true.

We then need to compute the nonlocal charge $z(P_2, v^n)$ defined in \cite{SZ} by
\begin{equation}\label{zh-10-6}
z(P_2, v^n)=\frac12g^{\alpha\beta}\nabla_\alpha(v^n)\nabla_\beta(v^n)-c\,v^n.
\end{equation}
It is equal to
\[\frac12 g^{nn}-\frac{1-d}2\,v^n=0.\]
This fact together with \eqref{zh-10-1} implies the locality of $\Phi(P_2)$.

Next, we need to show the coincidence of the metrics of $\Phi(P_i)$ and of $\hat{P}_i$, $i=1, 2$, this follows from the following identities
\begin{align*}
&(v^n)^{2}\hat{\eta}_{\alpha\beta}\,d\hat{v}^\alpha d\hat{v}^\beta=\eta_{\alpha\beta}\,dv^\alpha dv^\beta,\\
&(v^n)^{2}\hat{g}_{\alpha\beta}(\hat{v})\,d\hat{v}^\alpha d\hat{v}^\beta=g_{\alpha\beta}(v)\, dv^\alpha dv^\beta,
\end{align*}
and the transformation rule of the metrics of hydrodynamic Jacobi structures \cite{Fer, SZ}. The proposition is proved.
\end{prfn}

\section{Virasoro Constraints of the Tau Functions }\label{sec-5}

There is an important class of solutions  of the principal hierarchy \eqref{PH}
which can be obtained by solving the following system of equations \cite{D1, DZ2}:
\begin{equation}\label{EL}
\sum \tilde t^{\alpha,p}\frac{\p \theta_{\alpha,p}(v)}{\p v^\gamma}=0,\quad \gamma=1,\dots, n,
\end{equation}
where $\tilde t^{\alpha,p}=t^{\alpha,p}-c^{\alpha,p}$ and  $c^{\alpha,p}$ are  some constants which are assumed to be zero except for finitely many of them.
These constants are required to satisfy the genericity conditions that there exist constants $v^1_0,\dots, v^n_0$ such that
\[\sum c^{\alpha,p}\p_\gamma\theta_{\alpha,p}(v)|_{v=v_0}=0,\]
and the matrix
\[(A^\sigma_\gamma)=(\sum_{\alpha,p} c^{\alpha,p} \p^\sigma\p_\gamma \theta_{\alpha,p}(v_0))\]
is invertible. Here $\p_\gamma=\frac{\p}{\p v^\gamma}$ and $\p^\sigma=\sum_{\xi} \eta^{\sigma\xi} \p_\xi$.
One can obtain in this way a dense subset of the set of analytic monotonic solutions of the principal hierarchy \eqref{PH}, see Sec.\,3.6.4 of \cite{DZ2} for details.
The tau function for the solution $v^1(t),\dots, v^n(t)$ satisfying \eqref{EL} can be chosen
to be
\begin{equation}\label{logtau}
\log\tau(t)=\frac12\sum_{\alpha,\beta,p,q} \tilde t^{\alpha,p}\,\tilde t^{\beta,q}\, \Omega_{\alpha,p;\beta,q}(v(t)).
\end{equation}
The validity of the defining relation \eqref{tau-def} follows from  \eqref{EL} and the
identity \eqref{zh-22} for the functions $\Omega_{\alpha,p;\beta, q}$.

\begin{prp}\label{prp-v}
Let $v(t)=(v^1(t),\dots, v^n(t))$ be a solution of the principal hierarchy \eqref{PH} given by \eqref{EL}, and
$\hat{v}(\hat{t})=(\hat{v}^1(\hat{t}),\dots, \hat{v}^n(\hat{t}))$ be the  solution of the principal hierarchy \eqref{NPH} defined via \eqref{IS},
\eqref{zh-24-1}, \eqref{zh-24-2}. Then $\hat{v}(\hat{t})$ satisfies the equations
\begin{equation}\label{NEL}
\sum \tilde{\hat{t}}^{\alpha,p}\frac{\p \hat{\theta}_{\alpha,p}(v)}{\p \hat{v}^\gamma}=0,\quad \gamma=1,\dots, n,
\end{equation}
where $\tilde{\hat t}^{\alpha,p}=\hat t^{\alpha,p}-\hat c^{\alpha,p}$ with
\[\hat c^{1,0}=0,\quad  \hat{c}^{1,p+1}=-c^{n,p},\quad \hat{c}^{i,p}=c^{i,p},\quad \hat{c}^{n,p}=c^{1,p+1}\]
for $i\ne 1,n,\ p\ge 0$, and $\hat{\theta}_{\alpha,p}$ are defined in \eqref{ntheta}. The associated tau function
\begin{equation}\label{nlogtau}
\log\hat{\tau}(\hat{t})=\frac12\sum_{\alpha,\beta,p,q} \tilde{\hat t}^{\alpha,p}\,\tilde{\hat t}^{\beta,q}\, \hat\Omega_{\alpha,p;\beta,q}(\hat v(\hat t))
\end{equation}
satisfies \eqref{legendre}.
\end{prp}

\begin{prf}
To prove the validity of \eqref{NEL} let us consider the case when $\gamma=n$, the proof for other cases is similar.
By using \eqref{IS}, \eqref{ntheta} and \eqref{zh-24-2} we have
\begin{align*}
&\sum \tilde{\hat{t}}^{\alpha,p}\frac{\p \hat{\theta}_{\alpha,p}(\hat{v})}{\p \hat{v}^n}\\
=&\sum \tilde{\hat{t}}^{\alpha,p} \left(-\frac12\sum_{\sigma\ne 1,n} v_\sigma v^\sigma \frac{\p}{\p v^1}+\sum_{\sigma\ne 1,n} v^\sigma v^n \frac{\p}{\p v^\sigma}
+(v^n)^2 \frac{\p}{\p v^n}\right)\hat{\theta}_{\alpha,p}\\
=&-\frac1{2 v^n}\sum_{\sigma\ne 1,n} v_\sigma v^\sigma \sum_{\alpha,p} \tilde{t}^{\alpha,p} \frac{\p\theta_{\alpha,p}}{\p v^1}+\sum_{\sigma\ne 1, n} v^\sigma
\sum_{\alpha,p} \tilde{t}^{\alpha,p} \frac{\p\theta_{\alpha,p}}{\p v^\sigma}\\
&\quad +v^n \sum_{\alpha,p}\tilde{t}^{\alpha,p}\frac{\p\theta_{\alpha,p}}{\p v^n}-\sum_{\alpha,p} \tilde{t}^{\alpha,p} \theta_{\alpha,p}+\hat{t}^{1,0}-\hat{c}^{1,0}\\
=&-\sum_{\alpha,p}\tilde{t}^{\alpha,p} \theta_{\alpha,p}+\hat{t}^{1,0}
=-\frac{\p\log\tau}{\p x}+\hat{x}=0.
\end{align*}
Here we used the relation \eqref{zh-20b} and the fact that
\begin{equation}\label{tau-a}
\frac{\p\log\tau(t)}{\p x}=\sum_{\alpha,p} (t^{\alpha,p}-c^{\alpha,p})\theta_{\alpha,p}(v(t))
\end{equation}
which follows from \eqref{zh-22}, \eqref{EL} and \eqref{logtau}.
The validity of the relation \eqref{legendre} follows from \eqref{zh-24-2}, \eqref{relation}, \eqref{logtau} and \eqref{tau-a}. The proposition is proved.
\end{prf}

In the case when the solution $F(v)$ of the WDVV equations is conformal, the tau function
\eqref{logtau} satisfies the Virasoro constraints \cite{DZ3, EX, LT}
\begin{equation}
\sum a_m^{\alpha,p;\beta,q} \frac{\p\log\tau}{\p t^{\alpha,p}}\frac{\p\log\tau}{\p  t^{\beta,q}}+
\sum b_{m;\alpha,p}^{\beta,q} \tilde t^{\alpha,p} \frac{\p\log\tau}{\p t^{\beta,q}}
+\sum c_{m;\alpha,p;\beta,q} \tilde t^{\alpha,p} \tilde t^{\beta,q}=0,\label{zh-24-1b}
\end{equation}
where $m\ge -1$, and the coefficients that appear in the above expressions are some constants determined by the monodromy data of the Frobenius manifold of $F(v)$,
they define a set of linear differential operators
\begin{align}
L_m=&\sum a_m^{\alpha,p;\beta,q} \frac{\p^2}{\p t^{\alpha,p}\p  t^{\beta,q}}+
\sum b_{m;\alpha,p}^{\beta,q}  t^{\alpha,p} \frac{\p}{\p  t^{\beta,q}}\nn\\
&+\sum c_{m;\alpha,p;\beta,q} t^{\alpha,p} t^{\beta,q}+\delta_{m,0}\, c
\end{align}
which give a representation of  the half branch of the Virasoro algebra
\begin{equation}
[L_i, L_j]=(i-j) L_{i+j}+n\frac{i^3-i}{12} \delta_{i+j,0},\quad i, j \ge -1.
\end{equation}
The first two Virasoro operators have the expressions
\begin{align}
&L_{-1}=\sum_{p\ge 1} t^{\alpha,p}\frac{\p}{\p t^{\alpha,p-1}}+\frac12\,\eta_{\alpha\beta}
t^{\alpha,0} t^{\beta,0},\nn\\
&L_0=\sum_{p\ge 0} (p+\frac12+\mu_\alpha) t^{\alpha,p} \frac{\p}{\p t^{\alpha,p}}+
\sum_{p\ge 1} \sum_{r=1}^p (R_r)^\beta_\alpha\, t^{\alpha,p} \frac{\p}{\p t^{\beta,p-r}}\nn\\
&\quad\quad
+\frac12 \sum_{p,q\ge 0} (-1)^q \,(R_{p+q+1})^\xi_\alpha\, \eta_{\xi\beta}\, t^{\alpha,p}\, t^{\beta,q}
+\frac14 \sum_{\alpha} \left(\frac14-\mu_\alpha^2\right).
\end{align}
See \cite{DZ3, DZ2} for the explicit expressions of $L_m, m\ge 2$. From Proposition \ref{prp-v} it follows that the tau function of the principal hierarchy \eqref{NPH} obtained from the tau function \eqref{logtau} of \eqref{PH}
via the action of the inversion symmetry of the WDVV equations  satisfies the Virasoro
constraints
\begin{equation}
\sum \hat a_m^{\alpha,p;\beta,q} \frac{\p\log\hat \tau}{\p \hat t^{\alpha,p}}\frac{\p\log\hat \tau}{\p \hat  t^{\beta,q}}+
\sum \hat b_{m;\alpha,p}^{\beta,q} \tilde {\hat t}^{\alpha,p} \frac{\p\log\hat\tau}{\p \hat t^{\beta,q}}
+\sum \hat c_{m;\alpha,p;\beta,q} \tilde {\hat t}^{\alpha,p} \tilde {\hat t}^{\beta,q}=0
\end{equation}
associated to the solution $\hat F(\hat v)$ of the WDVV equations.

\section{The Topological Deformations} \label{sec-6}

The principal hierarchy \eqref{PH} possesses the following Virasoro symmetries \cite{DZ3}
\begin{equation}
\frac{\p \tau}{\p s_m}=
\sum a_m^{\alpha,p;\beta,q} \frac1{\tau} \frac{\p\tau}{\p t^{\alpha,p}}\frac{\p\tau}{\p  t^{\beta,q}}+
\sum b_{m;\alpha,p}^{\beta,q} t^{\alpha,p} \frac{\p\tau}{\p t^{\beta,q}}
+\sum c_{m;\alpha,p;\beta,q} t^{\alpha,p} t^{\beta,q} \tau,\label{zh-10-2}
\end{equation}
where $m\ge -1$. Note that these symmetries are nonlinear in $\tau$.

It is proved in \cite{DZ2} that, for a calibrated semisimple\footnote{A solution of the WDVV equation is called semisimple, if for any point $v$ of the associated Frobenius manifold $M$, the associative algebra defined on $T_vM$ by the structure constants $c_{\alpha\beta}^\gamma(v)$ is semisimple. This notion is not explicitly used in the present paper.} conformal solution $F(v)$ of the WDVV equations, there exists a unique deformation of the principal hierarchy
such that
\begin{itemize}
\item The deformed hierarchy possesses tau functions;
\item The Virasoro symmetries \eqref{zh-10-2} is deformed to the following one
\[\frac{\p \tau}{\p s_m}=L_m\tau,\quad m=-1,0,1,\dots.\]
which is linear in $\tau$.
\end{itemize}

The tau function of the deformed hierarchy that is specified by the Virasoro constraints
\[\left.L_m\right|_{t^{\alpha,p}\to t^{\alpha,p}-\delta^\alpha_1\delta^p_1} \tau(t)=0,\quad
m=-1,0,1,\dots.\]
corresponds to the partition function of a 2D TFT
if  the solution $F(v)$ of the WDVV equations
is given by the primary free energy of the matter sector of the 2D TFT.
Due to this fact such a deformation of the principal hierarchy is called the topological deformation, it has the form
\begin{equation}\label{zh-10-3}
\frac{\p w^\alpha}{\p t^{\beta,q}}= \eta^{\alpha\gamma}
\p_x\left(\frac{\p\theta_{\beta,q+1}(w)}{\p w^\gamma}\right)+
\sum_{g\ge 1} \varepsilon^{2 g} K^\alpha_{\beta,q;g}(w;w_x,\dots,w^{(2g+1)}).
\end{equation}
Here $\alpha,\beta=1,\dots,n$, $K^\alpha_{\beta,q;g}$ are differential polynomials of $w^1,\dots, w^n$, i.e. polynomials
of the $x$-derivatives of $w^1, \dots, w^n$ whose coefficients depend smoothly on $w^1, \dots, w^n$.

In this section, we will study the
relation between the topological deformations of the principal hierarchies of calibrated conformal solutions of the WDVV equations that are related by the
inversion symmetry.

We redenote the Hamiltonian structures $P_i$, $\{\ ,\ \}_i$, $H_{\alpha,p}$, \dots of the
principal hierarchy \eqref{PH} that appear in the previous sections by $P_i^{[0]}$, $\{\ ,\ \}_i^{[0]}$, $H_{\alpha,p}^{[0]}$, \dots
respectively. Then the topological deformation \eqref{zh-10-3} of the principal hierarchy also possesses a Hamiltonian structure given by the following data:
\begin{itemize}
\item[i)] A Hamiltonian operator $P_1$ with leading term $P_1^{[0]}$
\[P_1=P_1^{[0]}+\e^2\,P_1^{[1]}+\e^4\,P_1^{[2]}+\cdots,\]
where $P_1^{[k]}\ (k\ge1)$ are matrix differential operators whose coefficients are differential polynomials of $w^1,\dots, w^n$ with $\deg P_1^{[k]}=2k+1$;
\item[ii)] A set of differential polynomials $\{h_{\alpha,p}(w,w_x,\dots)\}$ of the form
\[h_{\alpha,p}(w,w_x,\dots)=\theta_{\alpha,p}(w)+\e^2\,h^{[1]}_{\alpha,p}(w,w_x,w_{xx})+\e^4\,h^{[2]}_{\alpha,p}(w,w_x,\dots)+\cdots,\]
where $h^{[k]}_{\alpha,p}\ (k\ge1)$ are differential polynomials of degree $2k$.
\end{itemize}
They define respectively the deformed Poisson bracket
\[\{H_1, H_2\}_1=\int \frac{\delta H_1}{\delta w^\alpha}\,P^{\alpha\beta}_1\left(\frac{\delta H_2}{\delta w^\beta}\right)dx
=\{H_1, H_2\}_1^{[0]}+\e^2\,\{H_1, H_2\}_1^{[1]}+\cdots\]
and the deformed Hamiltonians
\[H_{\alpha,p}=\int h_{\alpha,p+1}\,dx=H^{[0]}_{\alpha,p}+\e^2\,H^{[1]}_{\alpha,p}+\cdots.\]
In particular, the densities $h_{\alpha,0}$ are given by (c.f. \eqref{zh-10-5})
\[ h_{\alpha,0}=\eta_{\alpha\gamma} w^\gamma,\quad \alpha,=1,\dots,n.\]
Then the deformed hierarchy \eqref{zh-10-2} has the expression
\[\frac{\p w^\beta}{\p t^{\alpha,p}}=\{w^\beta, H_{\alpha,p}\}_1=P_1^{\beta\gamma}\,\frac{\delta H_{\alpha, p}}{\delta w^\gamma},\quad \alpha, \beta=1,\dots,n,\ p\ge 0.\]
We also denote $\p_{\alpha,p}=\frac{\p}{\p t^{\alpha,p}}$.
This hierarchy has the following properties:
\begin{itemize}
\item[i)] $\p_{1,0}=\p_x$;
\item[ii)] $\{H_{\alpha,p},H_{\beta,q}\}_1=0$, \  \ $\{H_{\alpha,-1}, \cdot\}_1=0$, \ \ $\p_{\alpha,p}H_{\beta,q}=0$,\  \ $[\p_{\alpha,p}, \p_{\beta,q}]=0$;
\item[iii)] $\p_{\alpha,p} h_{\beta,q}=\p_{\beta,q} h_{\alpha,p}$;
\end{itemize}
The property iii), which is called the tau symmetry condition, implies that for any pair of indices $(\alpha,p), (\beta,q)$ there exists a
differential polynomial $\Omega_{\alpha,p;\beta,q}$ such that
\[\Omega_{\alpha,p;\beta,q}=\Omega_{\alpha,p;\beta,q}^{[0]}+\e^2\Omega_{\alpha,p;\beta,q}^{[1]}+\cdots,\quad \p_x\Omega_{\alpha,p;\beta,q}=\p_{\alpha,p} h_{\beta,q}.\]
They are related to the tau function of the topological deformation of the principal hierarchy by
\[\Omega_{\alpha,p;\beta,q}(w(t),w_x(t),\dots)=\e^2\p_{\alpha,p}\p_{\beta,q}\log\tau(t)\]
where $w^{\alpha}(t)=\e^2\,\eta^{\alpha\beta}\p_{1,0}\p_{\beta,0}\log\tau(t)$. It follows from the definition of tau function $\tau^{[0]}(t)$ for the principal hierarchy
that
\begin{equation}\label{zh-10-9}
\F(t)=\e^{-2}\F_0(t)+\F_1(t)+\e^2 \F_2(t)+\cdots,
\end{equation}
where the free energy
$
\F(t)=\log\tau(t)$, $\F_0(t)=\log\tau^{[0]}(t)$.

The topological deformation of the principal hierarchy is constructed in \cite{DZ2}
by using the fact that the free energy $\F(t)$ can be determined by the requirement of the
linearization of the Virasoro symmetries via the genus zero free energy $\F_0$. Namely,
it can be represented in the form
\begin{equation}\label{zh-10-8}
\F(t)=\log\tau=\e^{-2} \F_0(t)+\Delta F(v,v_x,\dots)|_{v^\alpha=v^\alpha(t)}.
\end{equation}
Here
\[\Delta F=F_1(v,v_x)+\e^2 F_2(v, \dots, v^{\mathrm{(4)}})+\cdots+\e^{2g-2} F_g(v, \dots, v^{(3g-2)})+\cdots\]
with the functions $F_g$ determined by the loop equation give in \cite{DZ2}, and
\[v^\alpha(t)=\eta^{\alpha\gamma} \frac{\p^2\F_0(t)}{\p x \p t^{\gamma,0}},\quad \alpha=1,\dots,n\]
satisfy the principal hierarchy. The topological deformation of the principal hierarchy
is then obtained by the following
coordinates transformation
\begin{equation}
w^\alpha=v^\alpha+\e^2 \eta^{\alpha\beta}\p_{1,0}\p_{\beta,0}\Delta F. \label{quasi-miura}
\end{equation}
This transformations is a particular quasi-Miura transformation. In general, a quasi-Miura transformation will transform objects (like Hamiltonians,
vector fields, Hamiltonian structures, \dots) with differential polynomial coefficients to objects with coefficients being rational functions of the jet variables.
In \cite{DZ2} the above mentioned Hamiltonian structure of the deformed hierarchy
is obtained from the first Hamiltonian structure of the principal hierarchy via the
quasi-Miura transformation \eqref{quasi-miura}. A proof of the polynomiality of the
topological deformation of the principal hierarchy and of the  Hamiltonian operator $P_1$ are given in \cite{shadrin} and \cite{shadrin-new}. For the second Hamiltonian structure of the principal hierarchy, the following conjecture is given in \cite{DZ2}.
\begin{cnj}[\cite{DZ2}] \label{conj-1}
The quasi-Miura transformation \eqref{quasi-miura} transforms the second
Hamiltonian structure $P_2^{[0]}$ of the principal hierarchy to a Hamiltonian structure $P_2$ with
differential polynomial coefficients.
\end{cnj}
We assume its validity in the remaining part of the present section, and call the bihamiltonian structure $(P_1, P_2)$  that is induced from $(P_1^{[0]}, P_2^{[0]})$
via the quasi-Miura transformation \eqref{quasi-miura} the topological deformation of the bihamiltonian structure $(P_1^{[0]}, P_2^{[0]})$.

Let $(F(v), \{\theta_{\alpha,p}(v)\})$ and $(\hat{F}(\hat{v}), \{\hat{\theta}_{\alpha,p}(\hat{v})\})$ be two fixed calibrated conformal solutions of the WDVV equations
which are related by the inversion symmetry, and denote the topological deformations of the corresponding principal hierarchies by
\[\p_{\alpha,p}w^{\beta}=P_1^{\alpha\beta}\left(\frac{\delta H_{\alpha,p}}{\delta w^\beta}\right),\quad
\hat{\p}_{\alpha,p}\hat{w}^{\beta}=\hat{P}_1^{\alpha\beta}\left(\frac{\delta \hat{H}_{\alpha,p}}{\delta \hat{w}^\beta}\right)\]
respectively. Introduce the following reciprocal transformation $\Phi$:
\begin{align}
&d\tilde{x}=\sum_{\alpha=1}^n\sum_{p\ge0} h_{\alpha,p} dt^{\alpha,p},\label{dfm-rt-1}\\
\begin{split}\label{dfm-rt-2}
&\tilde{t}^{1,0}=\tilde{x},\quad \tilde{t}^{1, p}=-t^{n,p-1},\quad p \ge 1,\\
&\tilde{t}^{n,p}= t^{1,p+1}\ (p \ge 0),\quad \tilde{t}^{i,p}=t^{i,p},\quad  2\le i\le n,\ p \ge 0,
\end{split}
\end{align}
and the new coordinates $\tilde{w}^\alpha=\eta^{\alpha\beta}\tilde{h}_{\beta,0}$, where
\begin{align}
&\tilde{h}_{1,0}=-\frac{1}{h_{1,0}},\quad \tilde{h}_{1,p}=-\frac{h_{n,p-1}}{h_{1,0}},\quad  p\ge1,\nn \\
&\tilde{h}_{\alpha,p}=\frac{h_{\alpha,p}}{h_{1,0}},\quad  2\le \alpha\le n-1,\ p\ge 0, \label{dfm-h}\\
&\tilde{h}_{n,p}=\frac{h_{1,p+1}}{h_{1,0}},\quad  p\ge 0.\nn
\end{align}
Here the densities $h_{\alpha,p}=h_{\alpha,p}(w,w_x, \dots)$ are chosen as \cite{DZ2}
\[h_{\alpha,p}(w,w_x, \dots)=\theta_{\alpha,p+1}(v)+\frac{\p^2 \Delta F(v,v_x,\dots)}{\p x\p t^{\alpha,p+1}}.\]
The variables $v, v_x, \dots$ that appear on the right hand side of the above equation can be represented in terms of the variables $w, w_x, \dots$ by using the quasi-Miura transformation \eqref{quasi-miura}.

Denote the components of $\Phi(P_1)$ in the coordinates system $(\tilde{w}^1, \dots, \tilde{w}^n)$ by
$\tilde{P}_1^{\alpha\beta}$, then we have (see \cite{SZ})
\begin{equation}
\tilde{\p}_{\alpha,p}\tilde{w}^{\beta}=\tilde{P}_1^{\alpha\beta}\left(\frac{\delta \tilde{H}_{\alpha,p}}{\delta \tilde{w}^\beta}\right),\quad \tilde{\p}_{\alpha,p}=\frac{\p}{\p\tilde{t}^{\alpha,p}}. \label{hier-tran}
\end{equation}
It is easy to see that $\tilde{P}_1,\ \{\tilde{H}_{\alpha,p}\}$ have the same leading terms with $\hat{P}_1,\ \{\hat{H}_{\alpha,p}\}$ respectively but their deformed parts are different.

\begin{thm}\label{thm62}
There exists a Miura type transformation such that the hierarchy \eqref{hier-tran} is transformed to the topological deformation $\{\hat{\p}_{\alpha,p}\hat{w}^{\beta}\}$ of the principal hierarchy
 for $\hat{F}(\hat{v})$.
\end{thm}

\begin{lem}
Under the reciprocal transformation \eqref{dfm-rt-1}, \eqref{dfm-rt-2}
the bihamiltonian structure $(P_1, P_2)$ of the topological deformation of the principal hierarchy is transformed to a local bihamiltonian structure of
\eqref{hier-tran}.
\end{lem}
\begin{prf}
We are to use Theorem \ref{locality} again. Denote by $\Lambda=\int h_{1,0}\,dx$, we need to show that
\[[P_i, \Lambda]=0,\quad  z(P_i,h_{1,0})=0.\]
Here the function $z$ is defined as in \eqref{zh-10-6}.
The first equality is a consequence of the quasi-triviality of the bihamiltonian structure $(P_1, P_2)$ (see \cite{DLZ}), since we have proved
$[P_i^{[0]}, \Lambda^{[0]}]=0$. The second equality is verified in the proof of Proposition \ref{prp-reci}. The lemma is proved.
\end{prf}

\vskip 1em

\begin{prfn}{Theorem \ref{thm62}}
From the above lemma it follows that both the hierarchy \eqref{hier-tran} and the topological deformation $\{\hat{\p}_{\alpha,p}\hat{w}^{\beta}\}$ of the principal hierarchy for $\hat{F}(\hat{v})$ possess local bihamiltonian structures, these bihamiltonian structures have the same leading terms which form a  semisimple\footnote{A bihamiltonian structure of hydrodynamic type is semisimple if the eigenvalues
of the $(1,1)$ tensor $r^i_j=(g_2)^{ik}(g_1)_{kj}$ are non-constant and distinct, where $g_1, g_2$ are
the metrics associated to the hydrodynamic bihamiltonian structures. The bihamiltonian structures
obtained from a semisimple conformal solution of the WDVV equation is always semisimple. This
notion is not explicitly used in the present paper either.} bihamiltonian structure
of hydrodynamic type.

From Theorem 2.5.7 of \cite{SZ} it follows that Miura type transformations
preserve Schouten-Nijenhuis  brackets, so if two bihamiltonian structures are related
by a Miura type transformations, then this transformation transforms a bihamiltonian vector field to a
bihamiltonian vector field. By using this fact and the result of Corollary 1.9 of \cite{DLZ}  we know
that in order to prove the equivalence of two hierarchies under Miura type transformations, we only
need to show that their bihamiltonian structures are equivalent. According to the general results of
\cite{DLZ, LZ}, two bihamiltonian structures of the type considered here with same leading terms
are equivalent if and only if their central invariants\footnote{The notion of central invariants is
introduced in \cite{DLZ, LZ} to characterize infinitesimal deformations of a bihamiltonian
structure of hydrodynamic type. They are defined by certain tensor coefficients that appear in the
first six terms of the deformed bihamiltonian structure. Since their explicit definition will occupy
too much space and we never use it in the present paper, we omit it. In fact, all deformations
considered in this paper have central invariants $\frac1{24}$.} coincide. It is proved in \cite{DLZ2}
that the topological deformation of the bihamiltonian structure of a principal hierarchy has central
invariants $\frac1{24}$, so we have
\[c_i(P_1, P_2)=\frac1{24},\quad  c_i(\hat{P}_1, \hat{P}_2)=\frac1{24},\ i=1, \dots, n.\]
On the other hand, it is  shown in \cite{SZ} that if a reciprocal transformation transforms a \textit{local bi-Jacobi} (i.e. bihamiltonian)
structure to a local one, then it preserves the central invariants, which implies
\[c_i(P_1, P_2)=c_i(\tilde{P}_1, \tilde{P}_2),\quad  i=1, \dots, n,\]
so $c_i(\tilde{P}_1, \tilde{P}_2)=c_i(\hat{P}_1, \hat{P}_2)$, $i=1, \dots, n$. The theorem is proved.
\end{prfn}

The above theorem only ensures the existence of the Miura type transformation relating the two integrable hierarchies. We now consider the explicit form of this transformation.

Note that the reciprocal transformation \eqref{dfm-rt-1}-\eqref{dfm-h} is defined in the same way as in the dispersionless case, so we have the following proposition.
\begin{prp}
Let $\tau(t)$ be a tau function of the topological deformation \eqref{zh-10-3}
of the principal hierarchy associated to a calibration $\{\theta_{\alpha,p}(v)\}$ of $F(v)$. Define
\begin{equation}\label{legendre-new}
\log \tilde{\tau}(\tilde{t})=\log \tau(t)-x\frac{\p \log \tau(t)}{\p x},
\end{equation}
then $\tilde{\tau}(\tilde{t})$ is a tau function of the hierarchy \eqref{hier-tran}. It satisfies
\[\tilde{\Omega}_{\alpha,p;\beta,q}(\tilde{w}(\tilde{t}))=\e^2 \tilde{\p}_{\alpha,p}\tilde{\p}_{\beta,q}\log \tilde{\tau}(\tilde{t}),\]
where $\tilde{w}^\alpha(\tilde{t})=\eta^{\alpha\beta}\tilde{\p}_{1,0}\tilde{\p}_{\beta,0}\log \tilde{\tau}(\tilde{t})$, and $\tilde{\Omega}_{\alpha,p;\beta,q}(w)$ are defined
as in \eqref{relation} with $\theta_{\alpha,p}(v)$ replaced by $h_{\alpha,p}(w)$.
\end{prp}

We denote $\tilde\F(\tilde t)=\log\tilde\tau(\tilde t)$, and
\[\tilde\F(\tilde t)=\e^{-2}\tilde\F_0(\tilde t)+\tilde\F_1(\tilde t)+\e^2 \tilde\F_2(\tilde t)+\cdots.\]
It is easy to see that $\tilde\F_0=\log \hat{\tau}^{[0]}$, which is the tau function of the principal hierarchy
$\{\hat{\p}_{\alpha,p}\hat{v}^{\beta}\}$ given in Proposition \ref{prp-b}. In what follows, we are to show
that the functions $\tilde{\F}_g\ (g\ge1)$ can be obtained from $\F_g\ (g\ge1)$, and they are
in fact differential rational functions\footnote{A differential rational function of
$\hat{v}^1,\dots, \hat{v}^n$ is a rational function in $\hat{v}^i_x, \hat{v}^i_{xx}, \dots$ with
coefficients being smooth functions of $\hat{v}^1,\dots, \hat{v}^n$.} of $\hat{v}^1, \dots, \hat{v}^n$.

We regard the two sides of \eqref{legendre-new} as power series in $\e$ with coefficients being functions of $t$,
then compare their Laurent coefficients. By using \eqref{zh-10-8} we know that the reciprocal transformation
\eqref{dfm-rt-1} can be represented in terms of tau function by
\begin{equation}
\tilde x=\e^2 \frac{\p\log\tau(t)}{\p x}=\frac{\p\F_{0}(t)}{\p x}+\e^2 \frac{\p \Delta \F(t)}{\p x}=\hat x+
\e^2\frac{\p \Delta \F(t)}{\p x},
\end{equation}
where $\Delta\F(t)=\sum_{g\ge 1} \e^{2g-2} \F_g(t)$, so we have
\begin{equation}
\tilde{\F}_g(\tilde{t})=e^{\e^2 \frac{\p \Delta \F(t)}{\p x}\frac{\p}{\p \hat{x}}}\tilde{\F}_g(\hat{t})
=\sum_{k\ge0}\frac1{k!}\left(\e^2 \frac{\p \Delta \F(t)}{\p x}\right)^k\frac{\p^k \tilde{\F}_g(\hat{t})}{\p \hat{x}^k}.
\label{Fg-expand}
\end{equation}
In particular, by using the Legendre transformations \eqref{legendre}, \eqref{legendre-new} and the fact that $\tilde{\F}_0(\hat{t})=\hat{\F}_0(\hat{t})$ we have
\begin{align*}
&\e^{-2}\tilde{\F}_0(\tilde{t})=\e^{-2}e^{\e^2 \frac{\p \Delta \F(t)}{\p x}\frac{\p}{\p \hat{x}}}\tilde{\F}_0(\hat{t})\\
=&\e^{-2}\tilde{\F}_0(\hat{t})+\frac{\p \Delta \F(t)}{\p x}\frac{\p \tilde{\F}_0(\hat{t})}{\p \hat{x}}
+\e^{-2}\sum_{k\ge2}\frac1{k!}\left(\e^2 \frac{\p \Delta \F(t)}{\p x}\right)^k\frac{\p^k \tilde{\F}_0(\hat{t})}{\p \hat{x}^k}\\
=&\e^{-2}\left(\F_0(t)-x \frac{\p\F_0(t)}{\p x}\right)-x \frac{\p \Delta \F(t)}{\p x}
+\e^{-2}\sum_{k\ge2}\frac1{k!}\left(\e^2 \frac{\p \Delta \F(t)}{\p x}\right)^k\frac{\p^k \tilde{\F}_0(\hat{t})}{\p \hat{x}^k}\\
=&\e^{-2}\F_0(t)-x \frac{\p\F(t)}{\p x}
+\e^{-2}\sum_{k\ge2}\frac1{k!}\left(\e^2 \frac{\p \Delta \F(t)}{\p x}\right)^k\frac{\p^k \tilde{\F}_0(\hat{t})}{\p \hat{x}^k}\\
=&\log \tilde{\tau}(\tilde{t})-\Delta\F(t)
+\e^{-2}\sum_{k\ge2}\frac1{k!}\left(\e^2 \frac{\p \Delta \F(t)}{\p x}\right)^k\frac{\p^k \tilde{\F}_0(\hat{t})}{\p \hat{x}^k},
\end{align*}
by using \eqref{Fg-expand} again, we obtain
\begin{align}
&\sum_{g\ge1}\e^{2g-2}\sum_{k\ge0}\frac1{k!}\left(\e^2 \frac{\p \Delta \F(t)}{\p x}\right)^k\frac{\p^k \tilde{\F}_g(\hat{t})}{\p \hat{x}^k}\nn \\
&\qquad =\Delta\F(t)-\e^{-2}\sum_{k\ge2}\frac1{k!}\left(\e^2 \frac{\p \Delta \F(t)}{\p x}\right)^k\frac{\p^k \tilde{\F}_0(\hat{t})}{\p \hat{x}^k}.
\label{recF}
\end{align}
By comparing the coefficients of powers of $\e$ on both sides of the above equation we can represent
$\tilde{\F}_g$ in terms of ${\F}_0, \F_1, \dots,\F_g$ as follows: 
\begin{align}
\tilde{\F}_1(\hat{t})=&\F_1(t),\label{zh-10-11}\\
\tilde{\F}_2(\hat{t})=&\F_2(t)-\frac{\p \tilde{\F}_1(\hat{t})}{\p \hat{x}}\frac{\p \F_1(t)}{\p x}
-\frac12\frac{\p^2 \tilde{\F}_0(\hat{t})}{\p \hat{x}^2}\left(\frac{\p \F_1(t)}{\p x}\right)^2 \nn\\
=&\F_2(t)-\frac1{2v^n}\left(\frac{\p \F_1(t)}{\p x}\right)^2,\ \dots. \label{zh-10-12}
\end{align}
Here we used the relations
\[\frac{\p\tilde{\F}_1(\hat{t})}{\p \hat x}
= \frac{\p\F_1(t)}{\p x}\frac{\p x}{\p\hat x}=- \frac{\p\F_1(t)}{\p x}
\frac{\p^2\tilde\F_0(\hat t)}{\p \hat x^2}\]
and
\[ \frac{\p^2\tilde\F_0(\hat t)}{\p \hat x^2}={\hat v}^n=-\frac1{v^n}.\]
Note that the summation on the left hand side of \eqref{recF} starts from $k=2$, so we can use the fact that
\[\frac{\p^k \tilde{\F}_0(\hat{t})}{\p \hat{x}^k}=\frac{\p^{k-2} \hat{v}^n(\hat{t})}{\p \hat{x}^{k-2}},\]
then every $\tilde{\F}_g(\hat{t})$ is a differential polynomial of $v^n, \F_1, \dots, \F_{g}$. Note that
$\F_1, \dots, \F_{g}$ are differential rational function of $v^1, \dots, v^n$, and $\hat{v}^1, \dots, \hat{v}^n$ are rational functions of $v^1, \dots, v^n$, so we see that $\tilde{\F}_g$ is a differential rational 
function of $\hat{v}^1, \dots, \hat{v}^n$.

Now let us denote by
\[\log{\hat\tau}(\hat t)=\hat\F(\hat t)=\e^{-2} \hat\F_0(\hat t)+\sum_{g\ge 0} \e^{2g-2} \hat\F_g(\hat t)\]
the free energy of the topological deformation of the principal hierarchy associated to
$\hat F(\hat v)$.
Then as in \eqref{zh-10-8} we can represent $\hat \F _g(\hat t)$ in the form
\[\hat\F_g(\hat t)=\hat F_g(\hat v(\hat t),\dots,\p_{\hat x}^{3g-2} \hat v(\hat t)),\]
where $\hat{v}(\hat{t})=(\hat v^1(\hat t),\dots \hat v^n(\hat t))$ is a solution of the principal hierarchy of $\hat{F}(\hat{v})$.

Now let us compare $\hat{\F}_g$ and $\tilde{\F}_g\ (g\ge1)$.
When $g=1$, from \eqref{zh-10-11} and the expression of the genus one free energy \cite{DZ} it follows that
\[\tilde\F_1(\hat t)
=\F_1(t)=\frac{1}{24} \det\left(c_{\alpha\beta\gamma}(v(t)) v^\gamma_x(t)\right)+G(v(t))).\]
Here $G(v)$ is the $G$-function associated to $F(v)$, see \cite{DZ}. By using the following relation between the $G$-functions for $F(v)$ and its inversion $\hat F(\hat v)$ given in \cite{STR}:
\[
\hat{G}(\hat{v})=G(v)+(\frac{n}{24}-\frac{1}{2})\log v^n
\]
and the identity
\[
\frac{1}{24}\log \det(\hat{c}_{\alpha\beta\gamma}(\hat{v})\hat{v}_{\hat{x}_0}^\gamma)
  =\frac{1}{24}\log \det(c_{\alpha\beta\gamma}(v)v_{x}^{\gamma})-\frac{n}{24}\log v^n,
\]
we have
\[\tilde\F_1(\hat t)=\hat \F_1(\hat t)-\frac12 \log{\hat v^n(\hat t)}+\frac12 \log(-1).\]
For higher genera, we present the following conjecture.
\begin{cnj}
The difference ${\G}=\sum\limits_{g\ge 1}\e^{2g-2} \left(\tilde\F_g-\hat\F_g\right)$
can be represented as
\[\G(\hat{w})=\G_1(\hat{w}^n)+\e^2 \G_2(\hat{w}^n)+\e^4 \G_3(\hat{w}^n)+\cdots,\]
where $\G_g(\hat w^n)$ are differential polynomials of $\hat w^n$.
Moreover, the differential polynomials $\G_g$ do not dependent on the particular solution $F(v)$ of the WDVV equations. \end{cnj}

We have shown above that the conjecture holds true at the genus one approximation with
\[\G_1(\hat{w}^n)=-\frac12\log\hat{w}^n+\frac12\log(-1).\]
At the
genus two approximation, we have verified the validity of  the conjecture for solutions of the WDVV equations that are associated to the Coxeter groups of type $I_2(k)$ and $A_3$ with
\[\G_2(\hat{w}^n)=\frac{\hat{w}^n_{\hat{x}\hat{x}}}{8 (\hat{w}^n)^2}
-\frac{(\hat{w}^n_{\hat{x}})^2}{12 (\hat{w}^n)^3}.\]
For the higher genera corrections $\G_g \ (g\ge 3)$ we do  not know their explicit expressions at the moment. It is interesting to give an
interpretation for the expressions of these functions.

Under the assumption of validity of the above conjecture, the Miura type transformation
between the hierarchy \eqref{hier-tran} and the topological deformation
$\{\frac{\p\hat{w}^{\alpha}}{\p \hat t^{\beta,q}}\}$ of the principal hierarchy
for $\hat{F}(\hat{v})$ is given by
\[\tilde{w}^\alpha=\hat{w}^\alpha+\eta^{\alpha\beta}\p_{\hat{x}}\p_{\hat{t}^{\beta,0}}
\left[-\frac{\e^2}2 \log\hat{w}^n+\e^4 \left(\frac{\hat{w}^n_{\hat{x}\hat{x}}}{8 (\hat{w}^n)^2}
-\frac{(\hat{w}^n_{\hat{x}})^2}{12 (\hat{w}^n)^3}\right)\right]+\dots.\]
Here $\eta^{\alpha\beta}\p_{\hat{t}^{\beta,0}} \hat{w}^n=\hat{w}^\alpha_{\hat{x}}$.
We note that after the above Miura type transformation the flows
$\frac{\p\hat{w}^{\alpha}}{\p \hat t^{\beta,q}}$ are transformed to the
evolutionary PDEs
\[\frac{\p \tilde w^\alpha}{\p \hat t^{\beta,q}}=K^\alpha_{\beta,q}(\tilde w,\tilde w_{\hat x}, \tilde w_{\hat x\hat x},\dots),\quad \alpha,\beta=1,\dots. n,\ q\ge 0,\]
then the hierarchy \eqref{hier-tran} is obtained by redenoting the spatial variable $\hat x$
and the time variables $\hat t^{\beta,q}$ by $\tilde x$ and $\tilde t^{\beta,q}$ respectively.

\vskip 0.5truecm
\noindent{\bf Acknowledgments.}
The authors thank Boris Dubrovin for helpful discussions and comments.
The work is partially supported by the National Basic Research Program of China (973 Program) No.2007CB814800,
the NSFC No. 10801084, 11071135 and 11171176, and the Marie Curie IRSES project RIMMP.

\end{document}